\documentclass[12pt,a4paper,reqno]{amsart}
\usepackage{cite}
\usepackage{caption}
\usepackage{subcaption}
\usepackage{amsmath,amssymb,amsfonts,amsthm,epsfig,graphicx,xcolor}
\usepackage{esint}
\usepackage{url}
\usepackage{mathrsfs}

\usepackage[margin=0.8in]{geometry}

\newtheorem{theorem}{Theorem}[section]

\newtheorem{remark}[theorem]{Remark}

\theoremstyle{definition}

\numberwithin{equation}{section}
\numberwithin{figure}{section}

\renewcommand{\le}{\leqslant}

\def\R{\mathbb R}
\def\N{\mathbb N}

\title[A simple mathematical model inspired by the Purkinje cells]{A 
simple mathematical model
\\ inspired by the Purkinje cells:\\
from delayed travelling waves\\
to fractional diffusion}

\author{Serena Dipierro}
\address[Serena Dipierro]{Dipartimento di Matematica, Universit\`a degli studi di Milano,
Via Saldini 50, 20133 Milan, Italy}
\email{serena.dipierro@unimi.it}

\author{Enrico Valdinoci}
\address[Enrico Valdinoci]{School of Mathematics and Statistics,
University of Melbourne,
813 Swanston Street, Parkville VIC 3010, Australia, 
and Istituto di Matematica Applicata e Tecnologie Informatiche,
Consiglio Nazionale delle Ricerche,
Via Ferrata 1, 27100 Pavia, Italy,
and Dipartimento di Matematica, Universit\`a degli studi di Milano,
Via Saldini 50, 20133 Milan, Italy}
\email{enrico@mat.uniroma3.it}

\begin{document}

\begin{abstract}
Recently, several experiments have demonstrated the existence
of fractional diffusion in the neuronal transmission occurring
in the Purkinje cells, whose malfunctioning is known to be related 
to the
lack of voluntary coordination and the appearance of
tremors.

Also, a classical mathematical feature is that (fractional)
parabolic equations possess smoothing effects, in contrast with
the case of hyperbolic equations, which typically exhibit shocks and discontinuities.

In this paper, we show how a simple toy-model of a
highly ramified
structure, somehow inspired by that of the Purkinje cells,
may produce a fractional diffusion via the superposition
of travelling waves that solve a hyperbolic equation.

This could suggest that the 
high ramification of the Purkinje cells
might have provided an evolutionary advantage
of ``smoothing'' the transmission of signals and avoiding
shock propagations (at the price
of slowing a bit such transmission).
Although an experimental confirmation of the possibility
of such evolutionary advantage
goes well beyond the goals of this paper, we think that
it is intriguing, as a mathematical counterpart, to
consider
the time fractional diffusion 
as arising from the superposition of delayed travelling waves
in highly ramified transmission media.

The case of a travelling concave parabola with sufficiently
small curvature is explicitly computed.

The new link that we propose between time fractional diffusion
and hyperbolic equation also provides a novelty with respect to
the usual paradigm relating time
fractional diffusion with parabolic equations in the limit.

This paper is written in such a way as to be of interest
to both biologists and mathematician alike. In
order to accomplish this aim, both complete explanations of the objects considered
and detailed lists of references
are provided.
\end{abstract}


\maketitle

\section{Introduction}

Anomalous diffusion is becoming
increasingly more popular
to describe complex systems, in which the conventional diffusion described
by Brownian motion is inadequate. Among the different types of anomalous
diffusion, a special role is played by fractional diffusion, both in time and space variables.
Recent experiments
(see~\cite{ANAS, RR1, blanco, ERME01, ERME93, GGG, RR2, HHH, Magin,Fidel2, MARIN, MF,
MIRA, NACA, ROSS, Fidel, Fidel3, 44} and the references therein)
confirmed the evidence of fractional diffusions in many
systems of biological interest, though a complete understanding of all the
phenomena involved is still to be found, and many theoretical aspects of fractional
diffusion still needs to be further investigated.

This article is devoted to the derivation of a time fractional diffusion
equation from a ``basic building-block'', that is given by a set
of classical travelling waves. The model that we investigate is inspired by the Purkinje cells,
in which the complex ramification of the medium produces a selected
delay in the transmission of the waves. 

{F}rom the mathematical point of view, the superposition of these delays produces
a fractional operator. In this way,
a new link between
a time fractional diffusion equation (a parabolic equation) and
a classical wave equation (a hyperbolic one) will be introduced. 
In spite of the fact that these two types of equations present qualitatively very different properties,
the highly ramified structure of the transmission network allows us to reduce one equation
to the other.

{F}rom the biological point of view,
the results described in this article
may give a possible explanation of the causes of some malfunction of the Purkinje
cells.

\subsection{Aims and results}
The goal of this note is to provide {\em a mathematical derivation from
basic principle} of
the fractional diffusion equation
\begin{equation} \label{PDE}
D^s_t u = \kappa\, c^s\,L^{2-s}\,\partial^2_x u
\end{equation}
in view of some recent experiments on the
Purkinje cells (and, in general, on neural structures), which seem to exhibit this type of nonlocal diffusion.
In equation~\eqref{PDE}, the function~$u$ plays
the role of a diffusive substance (the concrete substance
depending on the particular type of diffusion considered in different
specific situations), $\kappa$ is an adimensional normalization
constant, and $c$ and~$L$ are fixed quantities (representing, respectively, the velocity of propagation of
an elementary signal and the length of the propagation device).

Moreover, the notation~$D^s_t$ stands for a time fractional derivative,
that, for definiteness, we take here in the sense of Caputo.
Namely, we recall the notion of Caputo fractional derivative
of order~$s\in(0,1)$ (see~\cite{caputo}), i.e. we set
\begin{equation}\label{EQ:C} 
D^s_t u(t):=
\frac{1}{\Gamma(1-s)} 
\int_0^t \frac{\dot u (\tau) }{(t-\tau)^{s}} \, d\tau,\end{equation}
where~$\Gamma$ is the Euler's Gamma-function
(which, for a fixed~$s\in(0,1)$, also plays in~\eqref{EQ:C}
just the role of a normalizing
factor).

Equations such as~\eqref{PDE} have recently appeared
in connection to several experimental data and theoretical considerations
related to the diffusion in the Purkinje cells: compare, in particular,
formula~\eqref{PDE} here
with the first formula in display in~\cite{Fidel2}
and see also~\cite{Fidel, EE}.

To the best of our knowledge, the scientific literature
has presented several deep and interesting descriptions
of the time fractional diffusion in~\eqref{PDE}
also in connection with neuronal biology, 
see e.g. the discussion related to formulas~(8.3), (8.12) and~(8.13)
in~\cite{MARIN}, but no attempt has been made till now to {\em derive
equation~\eqref{PDE} from ``basic principles''}
in a (possibly highly simplified)
toy-model somehow related to neurons.

Our goal in this paper is to try to fill this gap in the literature,
since we believe that derivations from basic principles
and from simpler equations have several
cultural and practical benefits, such as clarifying a difficult
but important research
subject, enlarging the
community of researchers working on a field,
providing connections with different subjects and
leading to a deeper and broader
understanding of the phenomena. Of course,
in this type of derivation processes some dramatic simplification
has sometimes to be expected, in order to reduce the arguments
to the core whenever possible, and, in this sense, the situation that
we present in this paper should not be intended as a ``full
explanation'' of the functioning of the complex neural networks,
but only as a simplified (though, in our opinion, sufficiently
``realistic'' and with concrete scientific value)
model, related to, or at least inspired by,
a simplified version of neural network.

Also, differently from the classical literature, our aim is not to relate
fractional diffusion with the standard one (which is formally obtained
in the limit as~$s\nearrow1$), but rather to
see fractional diffusion as a superposition of
hyperbolic equations with a delay.
In our setting, such delay is caused by the ramification
of the mathematical
structure on which the hyperbolic equation takes place.

The construction of this ramified medium is inspired by the structure of
the Purkinje cells,
which are a class of neurons located in the cerebellum, with
a highly ramified structure,
whose activity is of crucial importance
for the coordination of  complex motions.

As a matter of fact, the malfunctioning of
the Purkinje cells may lead (among other symptoms) to
ataxia (i.e., lack of voluntary coordination), 
tremors and hyperreactivity, see e.g.~\cite{blanco, MF}
and the references therein.

Thus,
one of the roles of the ramified structure of
the Purkinje cells seems to be that of somehow
``smoothing out''
sharp impulses. A mathematical counterpart of this phenomenon
can be seen by comparing the smoothing effects of the heat equation
with the shocks typical of hyperbolic equations (see e.g.~\cite{evans}).
Motivated by these considerations,
we provide a toy-model in which the fractional diffusion in~\eqref{PDE}
comes from the superposition of travelling waves.
In a sense, the highly ramified structure of the diffusion device (in the appropriate space/time
scale limit)
provides, at the end of the structure,
an averaged superposition of travelling waves with a nonlinear
delay which, in turn, transforms the hyperbolic equation of
a single travelling wave into a nonlocal {(in time)}
diffusion of the averaged
function as in~\eqref{PDE},
thus providing a regularizing effect on the solution.

This construction suggests the possibility
that the fractional diffusion experimentally
found in neurons could be related to the possibility of smoothing
irregular signals, so to make the coordination of the movements less subject to shocks and discontinuities.
In this sense, it is 
intriguing to wonder whether
the regularizing effects of fractional parabolic
equations (when compared to hyperbolic equations) could be seen as
a mathematical counterpart of
an evolutionary advantage of the high ramification of
the Purkinje cells, with the benefit of
smoothing the transmission of the signals and
perhaps favoring, at least indirectly,
a general coordination of the organism.

Roughly speaking, the mathematical effect of highly ramified arbors
in the transmission of signals may be thought as 
producing selected and appropriate delays
in the signal transmission.
The appropriately tuned superposition of
these delays has the combined effect to somehow ``average''
the propagation of
the signal, by producing a situation similar to those of fractional diffusion,
in which the speed at which diffusion takes place is ``anomalous'', i.e. it does not
coincide with the classical one prescribed by the Gaussian function.
Such
a delay behavior, when finely tuned, might somehow contribute to coordinate these signals, rather than just
retarding the whole transmission,
since the regularizing features of fractional parabolic equations
smooth out the data (differently from the case of standard
transmissions through hyperbolic equations).

{
\subsection{Fractional diffusion under different perspectives}

The contemporary literature in mathematical biology has considered
anomalous diffusion of fractional type in several experiments and under different
points of view.
Several recent works considered
fractional diffusion in view of
the so-called
``input-output analysis'': namely, accurate measurements
are performed to discover
underlying biological mechanisms, often related to linear equations.
See in particular~\cite{Thorson} for
the point of view of 
physiology on the power law distributions arising in the
outputs of several receptors.
In this, a classical approach of~\cite{Schweidler}
is that of considering powers of time as
linear superpositions of
many different exponential decays,
in conformity with the definition of the Euler's Gamma-function,
see in particular formula~(1)
in~\cite{Thorson} and the references therein (a Gamma-function approach
is also useful to link space fractional diffusion and classical heat semigroups, see
e.g. formulas~(2.11) and~(2.12) of~\cite{bucur2}).

The approach of~\cite{Schweidler, Thorson} has also been exploited
in the analysis of cells and tissues, see e.g. formula~(2.2) in~\cite{Magin}.
Also, in~\cite{Magin} fractional calculus is used
to study a tissue-electrode interface between cardiac muscle cells.
Related methods have also been exploited in~\cite{ANAS}
for the analysis of the vestibulo-ocular system.

In~\cite{DUHU}, the time fractional derivative is used to model several phenomena.
First, fractional derivatives of order~$s\in(0,1)$ are seen as natural interpolations
between elasticity laws of Hooke type (in which the displacement corresponds to
a derivative of order~$s=0$) and viscosity theories of Newton type (in which the velocity
corresponds to
a derivative of order~$s=1$): with this respect, fractional derivatives turn out to provide an
interesting model for viscoelastic materials. Then, some biological models
are discussed in~\cite{DUHU}, with special attention to protein adsorption kinetics.
Finally, cognitive processes are also taken into account in~\cite{DUHU}, showing
a fit with classic memorizing testing data by Hermann Ebbinghaus which date
back to~1885.
As a technical remark, we observe that in~\cite{DUHU}
the Riemann-Liouville derivative is taken as basic model for time fractional diffusion,
instead of the Caputo derivative that we consider here (compare formula~(2)
in~\cite{DUHU} with~\eqref{EQ:C} here): nevertheless the two fractional derivatives
are closely related, up to a term involving the initial condition, see e.g. formula~(4)
in~\cite{EQUI}.

A recent review
of many different results of
fractional calculus in bioscience and engineering, with special emphasis
in respiratory tissues and drug diffusion is given in~\cite{IONE}.
See in particular Section~4.3 in~\cite{IONE} for a discussion on
the many applications of fractional calculus in neuroscience.

Though our paper mainly focuses on one single example,
namely that of the transmission problem in a highly ramified network,
we believe that
our approach is general enough to be applicable also to other models (see e.g.
Section~4.4 in~\cite{ABAR}).
}

\subsection{Topics and methods}
In our analysis, the mathematical methods exploited are all of
elementary nature (integration by parts, change of variable,
superposition principle, integration theory, basic PDE), 
all the ansatz and approximation assumptions are clearly stated
step by step,
the biological motivations are explained without assuming
major prerequisites and providing quite exhaustive lists of classical
and contemporary references, and
we also indulge in explanations and clarifications,
therefore the
paper is {\em easily usable by a wide group of interested researchers}.
Our goal here is not to give a complete explanation
of the rich phenomena encoded by the complex structure of
the Purkinje cells; nevertheless, we believe that the mathematical
approach that we present here
is useful to better understand some specific neuronal
features. Moreover, the new connection between fractional diffusion
and hyperbolic equations may lead to a better understanding of
the smoothing effects for travelling signals that are favored by the highly
ramified neuronal structures.

Furthermore, though a clear understanding of the time fractional diffusion
in real world phenomena
will require combined efforts from different perspectives
(e.g. with synergic approaches from biology, chemistry, physics, etc.),
we hope that the mathematical insight presented here
can better motivate the interest of fractional diffusion among the
mathematical community, serve as a foundation
ground for scientists with different backgrounds and suggest
new connections between very different types of evolution equations,
such as hyperbolic and (fractional) parabolic ones, which are made possible
by the highly complex structures of the media. Also, we believe that a
mathematical model easy to handle
may provide some initial insight (to be enhanced
by quantitative and more sophisticated
studies)
about the level at which
fractional diffusion arises in many natural phenomena,
also trying to give information on the scales
involved and on the basic causes of these
features.

\subsection{Organization of the paper}

The rest of this paper is organized as follows.
In Section~\ref{KAJ} we recall some basic facts about the Caputo derivative
and the related time fractional diffusion. In particular,
a simple integration by parts procedure, combined
with an appropriate change of variables, relates the
Caputo derivative to the superpositions of delayed classical second
derivatives.

Then, in Section~\ref{tgsfd0202}, we will introduce a highly ramified
mathematical
structure, inspired by the neuronal spikes,
and study the transmission of a hyperbolic travelling wave
along such complex medium. We will see that the structure
of the medium produces the superposition of delayed travelling waves
which, at the end of the transmission
device, in average and in the appropriate limit sense, can be related to
the fractional derivative and produce the time fractional diffusion equation
in~\eqref{PDE}.

In the computations needed for this scope, an ancillary limit formula
is stated, whose proof is given, for the facility of the reader, in Section~\ref{STBD}.

The conclusions of this paper are then summarized in Section~\ref{STBD:2}.

\section{Integration by parts in the Caputo derivative}\label{KAJ}

The notion of fractional diffusion provides at the moment
an intense topic of research, both for its very challenging theoretical difficulties
and in view of concrete applications in biology, physics and finance
(see e.g.~\cite{bucur2} for several explicit discussions
and motivations). In general, fractional diffusion presents
several phenomena in common with the classical diffusion
arising from Gaussian processes and Brownian motions, such as
the regularizing effects with respect to initial data:
in this sense, see in particular~\cite{kim}
for a regularity theory in Lebesgue spaces,
and also~\cite{zacher} for a regularity
theory in H\"older spaces (see also
page~103 in~\cite{zacher2}
for a general discussion
about the relation between ``abstract
Volterra equations'' and the time fractional diffusion).
For higher regularity in time, see also
Section 5 of~\cite{allen}, and for related results see~\cite{allen-et}.

Furthermore, important and often unexpected differences
between classical and fractional diffusion arise:
for instance, solutions of fractional equations can locally approximate
any given function, in sharp contrast with the classical case,
see~\cite{JJ, bucur, JJ2}.

The monograph~\cite{kai} also provides extensive and throughout discussions
about fractional derivatives in time also in view of many applications.
See also~\cite{YY77} for an approach to different types
of time fractional derivatives from the
perspectives of stochastic processes with {long rests}.

Here, we recall some basic facts on the fractional
Caputo derivative in~\eqref{EQ:C}
and perform some preliminary computations
which will be used in the forthcoming sections.
To start with, for notational convenience, we scale the constant
in~\eqref{EQ:C} and we integrate by parts, by
obtaining that
\begin{equation}\label{EG:Q2}
\begin{split}
\partial^s_t u(t)\,&:= (2-s)\,(1-s)\,\Gamma(1-s)\,D^s_t u(t)
\\ &=(2-s)\,(1-s)\,
\int_0^t \frac{\dot u (\tau) }{(t-\tau)^{s}} \, d\tau\\
&= 
-(2-s)\,\int_0^t 
\frac{\partial}{\partial\tau} 
\big( (t-\tau)^{1-s}\,\dot u (\tau)  \big)\, d\tau
+(2-s)\,\int_0^t  (t-\tau)^{1-s}\,\ddot u (\tau) \, d\tau\\
&=(2-s)\,t^{1-s}\,\dot u(0)
+(2-s)\,\int_0^t  (t-\tau)^{1-s}\,\ddot u (\tau) \, d\tau.
\end{split}
\end{equation}
Using the substitution~$\vartheta:=(t-\tau)^{2-s}$, we thus obtain
\begin{equation}\label{EG:Q3}
\partial^s_t u(t)=(2-s)\,t^{1-s}\,\dot u(0)
+\int_0^{t^{2-s}}\ddot u (t-\vartheta^\beta) \, d\vartheta,
\end{equation}
where
\begin{equation}\label{EG:Q5} \beta:=\frac1{2-s}\in\left(\frac12,1\right).\end{equation}
{As a matter of fact,
the computations leading to~\eqref{EG:Q2} are merely formal, since we are
tacitly assuming here that~$u$ is smooth and has a well-defined second time derivative.}
Our goal is now to interpret~\eqref{EG:Q3}
as a superposition of delayed effects caused by the ramified
structure of the transmission medium, which is somehow inspired by the structure
of the Purkinje cells (in neuronal
transmissions, other types of delays leading to fractional diffusion
may be caused by obstacles and bindings, see~\cite{BIND}).

\section{A simple model towards the fractional diffusion
in the Purkinje cells}\label{tgsfd0202}

In this section, we present a transmission media built by a highly ramified
structure. Such mathematical model is qualitatively
inspired by the dendritic arbor
of the Purkinje
cells. We will consider the transmission along this medium,
as prescribed by the classical wave equation.
The ramifications of the medium will cause the delay of some signals,
whose superposition at the end of the structure will be related to the
fractional derivative.

In this setting, the superposition of travelling waves with 
a suitably tuned delay will produce a time fractional diffusion
in the formal limit,
thus providing a new bridge between equations of very different kind
(i.e. hyperbolic and fractional parabolic) with the aim of encoding
some of the features observed by the experiments in neuronal dendrites,
such as the time fractional diffusion in the spikes of the
Purkinje
cells.

In the classical
transmission line analysis (see e.g.  pages~5--14 in~\cite{KARA})
the variation in space of the potential
is related, via inductance,
to the variation in time of the intensity of current;
on the other hand,
the variation in space of the intensity of current
is related, via capacitance,
to the variation in time of the potential. The combination of these equations
naturally lead to the wave equation
(for a derivation of the wave equation directly
from Maxwell's equations see e.g.
(3.14) and~(3.15) in~\cite{COST}).
For models presenting wave equations in cylindrical neurons,
see e.g. formulas~(6.20) and~(6.21) in~\cite{RIGA},
and also~\cite{RR1, RR2, RR3, RR4}.

We remark that the analysis of travelling waves 
in neurons is a classical and active topic of study in itself,
see e.g.~\cite{ERME93, ERME01, ERME01b, MIRA, COO}.
See also~\cite{44} for other mathematical models
related to diffusion in neurons based on Monte Carlo methods.

\subsection{Model of the ramified structure}
The model of the ramified medium that we take into
account is a very basic simplification inspired by a branch of the Purkinje
cell and goes as follows.

We let~$N\in\N$ (to be taken large in the sequel)
and
\begin{equation}\label{8usdhjhdgdsdf13vh:0} 
b_N:=\sum_{1\le k\le N} \frac1{{k}^{\alpha}},\end{equation}
with
\begin{equation}\label{8usdhjhdgdsdf13vh:0X89} 
\alpha:=\frac{1-s}{2-s}=1-\beta\in \left(0,\frac12\right)\end{equation}
and~$\beta$ as in~\eqref{EG:Q5}. We let also
\begin{equation}\label{8usdhjhdgdsdf13vh}
\begin{split}& \ell_k:= \frac{L}{k^\alpha b_N}\\{\mbox{and }}\quad&
\lambda_k:=\sum_{1\le j\le k}\ell_j.\end{split}\end{equation}
We consider a set of~$N$ planar curves, with one common
endpoint and the other endpoint lying on a common straight line.
These curves will be denoted by~$S_1,\dots,S_N$.
The length of~$S_1$ is set to be~$L>0$.
The length of~$S_2$ is set to be~$L+\ell_1$.
Iteratively, the length of~$S_k$, for each~$k\in\{2,\dots,N\}$
is set to be equal to
$$ L+\ell_1+\dots+\ell_{k-1}=L+\lambda_{k-1}.$$

\begin{figure}
    \centering
    \includegraphics[height=6.9cm]{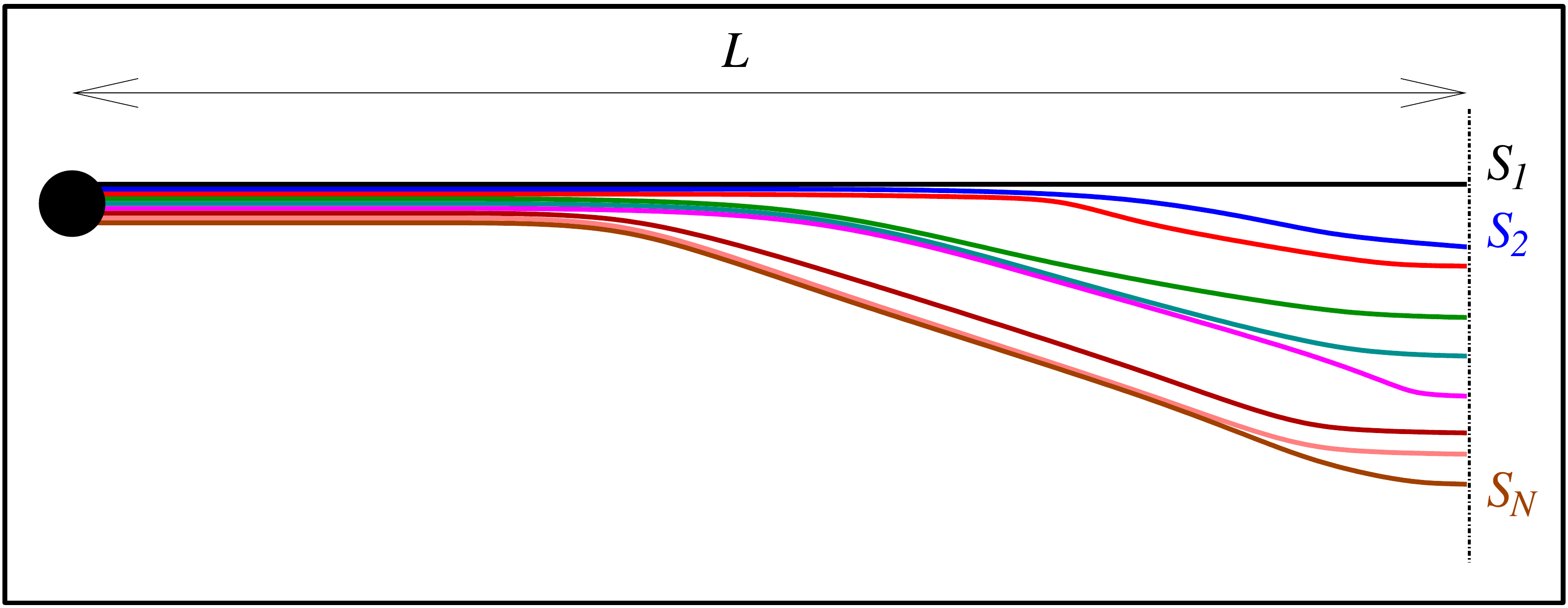}
    \caption{The ramified media along which 
we study the delayed transmission (an exemplifying,
not quantitative, picture).}
    \label{NESS}
\end{figure}

See Figure~\ref{NESS} for an exemplifying representation of
such medium.
Of course, the reader may compare this picture and the classical drawing
by S. Ram\'on y Cajal, see e.g.~\cite{WIKI}
as well as the many
realistic pictures available nowadays,
to appreciate the motivation related to the
Purkinje
cell -- though of course we do not aim to a full
understanding of the Purkinje cell by the very simplified
model described here. 
We stress, in particular, that when we refer to concepts like ``horizontal'',
``left'' or ``right'' we are clearly making no reference to concrete\footnote{As
a matter of fact, from the biology perspective,
it also makes sense of looking at signal travelling from the right
to the left in a medium as the one described in Figure~\ref{NESS}.
In this case, the left end of Figure~\ref{NESS}
would act as a ``soma'' (the cell body of a neuron),
which is the site of final integration of the
signals and ultimately is responsible to
generate action potentials to send to downstream neurons.
In our mathematical deduction, this
situation can be also taken into account
(what counts is just the average delayed obtained
by a signal travelling from one end of the transmission
medium to the other).}
special directions in the cerebellum, but rather we are referring
to the toy-model drawn in Figure~\ref{NESS}:
in any case, as it will be apparent in the computation,
these directions play no role in our derivation, what counts
is only the fact that the ramification produces branches of different
lengths or, better to say, just the fact that the speed to
travel to the end of a branch depends on the branch itself.
Also, Figure~\ref{NESS} can be immersed in a more complicated
neural network and we only focus on a ``single'' ramified system
to keep the discussion as simple as possible (in practice,
of course, the response to an organism depends on the
signal transmissions in many cells and not just one).

See also~\cite{fiala, 2002} and the references
therein for accurate descriptions of dendritic arbors.

Of course, different scalings and different parameters in the model
presented
lead to quantitatively different\footnote{Other models
with different structures
may also be taken into account
to address anomalous diffusion in related contexts.
For instance, 
it is interesting to
investigate whether
some sort of hyperbolic superposition with delay 
can also be applied to cases that do not
need the presence of multiple branches. In particular,
in~\cite{Fidel3} the pyramidal 
cells in the hippocampus are studied,
showing that the presence of dendritic 
spines causes anomalous diffusion: indeed, 
these small protrusions force diffusing molecules 
to undergo a continuous random walk 
with random waiting times that result in anomalous diffusion.  }
versions of fractional diffusions,
which may be seen as the mathematical
counterpart of the quantitative differences explored in~\cite{Fidel3}
in view of the different densities of dendritic spines.
See also~\cite{ben} for an accurate analysis
of anomalous diffusion in
fractal media. Of course, a detailed description of the geometry
of the dendritic spines also in terms of density and statistics of
the ramification may produce a better understanding of the
neuronal diffusion and inspire quantitatively
more accurate models.
In addition, it would be very interesting to
further investigate these questions and related problems
in the light of the theory of wave propagation in networks
(see e.g.~\cite{zuazua} and the references therein,
where however Dirichlet or control conditions are usually assumed on the
vertices of the graphs, also in relation with vibrating
strings with joined extrema).

\subsection{Travelling waves in the ramified medium and their effect
at the right end of the structure}

We will now take into account a travelling wave
in the structure depicted in Figure~\ref{NESS},
and compute the averaged effect of the delay induced 
by the ramification of the medium.

{F}rom the biological point of view,
we think it is worth justifying our choice of considering
hyperbolic equations as ``building blocks''
of our derivation procedure.
Indeed, some of the equations proposed to model
neuronal diffusion, such as the Hodgkin-Huxley equation,
are of parabolic, rather then hyperbolic type.
Nevertheless, we chose not to consider the Hodgkin-Huxley equation
as the basis of our model, for several reasons.
First of all, the Hodgkin-Huxley equation~\cite{HHH}
is a rather complicated model and, at least from an ``aesthetic''
perspective, does not seem to be suitable as a ``basic principle''
from which one derives a time fractional process.
More importantly, 
there are in the literature
some non-thermodynamical 
theories 
that suggest (and some experiments support) 
that mechanical waves accompany action potentials, see~\cite{EL, GGG}.
Furthermore,
it is known that
some reaction-diffusion processes 
can produce solitons and travelling waves, see~\cite{SSS}.

For these reasons, we thought it was intriguing, and also
sufficiently realistic, to take the wave equation as the
basic to deduce a time fractional equation
in our simplified setting.

Thus, the idea is to consider a travelling wave 
(say, from left to right) in such medium,
whose ``natural speed'' is given by~$c>0$, 
that is a (smooth)
function
\begin{equation}\label{fo}
f=f(x,t)=f_o\left(
{
\frac{x-ct}{L} }\right)\end{equation} satisfying\footnote{For the sake of simplicity, here
we take into account equations in the simplest possible form: in particular,
we are not explicitly considering forcing terms in the equation (which can be anyway
added in a more general analysis).
Also, in order to develop formal expansions, we assume~$f_o$ to be smooth,
with a smoothness independent of the structural parameters~$N$, $L$ and~$c$.
{The advantage of the scaling in the variable of~$f_o$ is that its dependence on
the spatial variables is weighted by the natural scale of the problem, thus~$f_o$
is a function of ``adimensional'' coordinates.}}
\begin{equation}\label{ADD:con:0} \partial^2_t f(x,t)=c^2 \partial^2_x f(x,t)={\frac{c^2}{L^2}
f_o''\left({\frac{x-ct}{L} }\right)}, 
\end{equation}
and analyze (for large~$N$)
its effect on the right end of the structure, which we will denote by~$u$.

In our toy-model, the function~$u$ at the right end
of the structure represents, somehow, the relevant information
that the left end of the structure ``sends'' to the organism:
in this sense, we think it is natural to consider such information
at the right end of the structure as the superposition of the information
sent in each of the branches which connect the left end to the right end.
In view of this,
we will see that~$u$ is the superposition of a series of~$f$'s,
shifted by a nonlinear delay (a precise formula will be given in~\eqref{0oxc:A}
below).

Roughly speaking, the effect of any travelling
function~$f$ on the right end of the structure
may be seen as the superposition of the~$1/N$-contributions
of $f$ along each of the ramifications~$S_k$, for~$k\in\{1,\dots,N\}$.
Each of these contributions will be denoted by~$f_k$.
We drop the dependence on~$x$ for the sake of 
simplicity
and we assume that~$f_1(t)=\frac{f(t)}{N}$, namely the contributions
are ``equally spread'' on the ramifications.
Also,
we observe that the length travelled by the function~$f_2$
is equal to the one travelled by~$f_1$ (which is in turn given by~$L$)
plus the length of the additional quantity~$\ell_1$. This causes a phase delay of~$f_2$
with respect to~$f_1$ of size~$ c^{-1}\ell_1$. Hence
$$ f_2(t)=f_1(t-c^{-1}\ell_1).$$
Iteratively, for each~$k\in\{2,\dots,N\}$,
the length travelled by the function~$f_k$
is equal to the one travelled by~$f_{k-1}$
plus the length of the additional quantity~$\ell_{k-1}$. This causes a phase delay of~$f_k$
with respect to~$f_{k-1}$ of size~$ c^{-1}\ell_{k-1}$. Hence
$$ f_k(t)=f_{k-1}(t-c^{-1}\ell_{k-1}).$$
Accordingly
\begin{eqnarray*}&& f_k(t)=f_{k-1}(t-c^{-1}\ell_{k-1})
=f_{k-2}(t-c^{-1}\ell_{k-1}-c^{-1}\ell_{k-2})\\ &&\qquad=\dots=
f_{k-i}(t-c^{-1}\ell_{k-1}-c^{-1}\ell_{k-2}-\dots-c^{-1}\ell_{k-i})\\ &&\qquad=\dots=
f_{1}(t-c^{-1}\ell_{k-1}-c^{-1}\ell_{k-2}-\dots-c^{-1}\ell_{1})=
f_1(t-c^{-1}\lambda_{k-1})\\ &&\qquad=\frac{f(t-c^{-1}\lambda_{k-1})}{N}.\end{eqnarray*}
Hence, the total contribution of a travelling function~$f$ on the right side
of the structure is taken to be
\begin{equation}\label{00TBD} \sum_{1\le k\le N} f_k(t)= \frac{1}{N}
\sum_{1\le k\le N} f(t-c^{-1}\lambda_{k-1}).\end{equation}
We set
\begin{equation}\label{uaso:n}\begin{split}&
\eta_{k,N}:= L\, \left(\frac{k}{N}\right)^{1-\alpha}-\lambda_{k-1}\\{\mbox{and }}\;&
\epsilon_N:=
\sup_{k\in\{1,\dots,N\}} |\eta_{k,N}| \end{split}\end{equation}
and it can be shown that
\begin{equation}\label{TBD}
\lim_{N\to+\infty}\epsilon_N=0.
\end{equation}
Not to interrupt this calculation, we postpone the proof of~\eqref{TBD}
to Section~\ref{STBD}.

Now, in view of~\eqref{00TBD}
and~\eqref{TBD}, the
total contribution of a smooth travelling function~$f$ on the right side
of the structure becomes
\begin{equation}\label{RIED}
\begin{split}
&\frac{1}{N}
\sum_{1\le k\le N} f\left(t-c^{-1}L\,\left(\frac{k}{N}\right)^{1-\alpha}+c^{-1}\,\eta_{k,N}\right)
\\ =\;&
\frac{1}{N} 
\sum_{1\le k\le N} \left[
f\left(t-c^{-1}L\,\left(\frac{k}{N}\right)^{1-\alpha}\right)
+O(c^{-1}\,\epsilon_N)\right]\\
 =\;&
\frac{1}{N} 
\sum_{1\le k\le N} 
f\left(t-c^{-1}L\,\left(\frac{k}{N}\right)^{1-\alpha}\right)
+O(c^{-1}\,\epsilon_N).
\end{split}\end{equation}
We now recognize a Riemann sum, namely we have that
\begin{equation} \label{7uw3e}\begin{split}&
\lim_{N\to+\infty}
\frac{1}{N} 
\sum_{1\le k\le N} 
f\left(t-c^{-1}L\,\left(\frac{k}{N}\right)^{1-\alpha}\right)\\=\;&
\int_0^1 f\left(t-c^{-1}\,L\,\zeta^{1-\alpha}\right)\,d\zeta.\end{split}
\end{equation}
{ For clarity, we now distinguish
between the time variable~$t$ and the natural
time scale of the problem, which will be denoted by~$T$. 
While~$ t$ is a free coordinate and it is one of the arguments 
of the functions under considerations, $T$ is the ratio between the 
characteristic length of the propagation device and the natural velocity 
of propagation of the signal, namely
\begin{equation}\label{Speed}
T:=\frac{L}{c}.\end{equation}
Notice that~$ T$ is fixed in terms of the medium and the propagation speed, 
hence can be considered as a characteristic feature of the system under consideration.
In this way, the functions can be differentiated in the variable~$t$,
and the result can be evaluated, for instance, at~$T$, with the aim of reconstructing 
a time fractional diffusion equation when the spatial scale is that of 
the end of the propagation device and the time scale is of the order of the characteristic time~$ T$.

We also define
\begin{equation}\label{DELTA}
\delta:=\frac{1}{T^s}=\left(\frac{c}{L}\right)^s.
\end{equation}
It is interesting to
point out that, if we wish,
we can
consider~$\delta$ as a small\footnote{{We observe that, with respect to the parameters~$c$ and~$L$
in~\eqref{PDE}, it is possible to choose~$\delta$ small, without making the equation
in~\eqref{PDE} degenerate. For instance, if one takes~$c\sim L^{\frac{s-2}{s}}$, then the coefficients in~\eqref{PDE} do not degenerate and~$\delta$ is small
if~$L$ is large.}}
parameter.

In a sense,
the setting in~\eqref{Speed} and~\eqref{DELTA}
says that equation~\eqref{PDE} can be obtained from
superposed hyperbolic equations ``only in the appropriate
space/time scaling''. Of course, we cannot reduce the mathematical
complexity of the hyperbolic equations; moreover,
this appropriate choice of
scaling might have a biological meaning in the neuron transmission,
since ``waves
are rarely detected beyond the point where the thick dendrites
begin to branch'', according to page~4 of~\cite{ROSS},
quoting~\cite{LARKUM, NACA}.

The setting in~\eqref{Speed} is exploited in combination with the substitution~$\vartheta:=
(c^{-1}\,L)^{\frac1{1-\alpha}}\zeta$. }
That is, recalling~\eqref{8usdhjhdgdsdf13vh:0X89}, we have that~$\vartheta={T}^{2-s}\zeta$
and
$$ \int_0^1 f\left(t-c^{-1}\,L\,\zeta^{1-\alpha}\right)\,d\zeta=
\frac1{{T}^{2-s}}\,
\int_0^{ {T}^{2-s}} f\left(t-\vartheta^\beta\right)\,d\vartheta.$$
Plugging this information into~\eqref{7uw3e} we obtain that
\begin{eqnarray*}&& \lim_{N\to+\infty}
\frac{1}{N} 
\sum_{1\le k\le N} 
f\left(t-c^{-1}L\,\left(\frac{k}{N}\right)^{1-\alpha}\right)\,=\,
\frac1{ {T}^{2-s} }\,
\int_0^{ {T}^{2-s}} f\left(t-\vartheta^\beta\right)\,d\vartheta.\end{eqnarray*}
{F}rom this, \eqref{EG:Q5}, \eqref{TBD},
\eqref{RIED} and~\eqref{Speed} we conclude that,
for large~$N$, we can approximate the
contribution of a smooth travelling function~$f$ on the right side
of the structure with the quantity
\begin{equation}\label{0oxc:A}
u(t):=\frac{1}{{T}^{2-s}}
\,
\int_0^{{T}^{2-s}} f\left(t-\vartheta^{\frac1{2-s}}\right)\,d\vartheta.
\end{equation}

\subsection{A fractional equation for $u$}
Now, we consider the case of propagation
in a bounded spatial region
and we will show that, in a suitable approximation, and at a space/time scale coherent
with that of the end of the transmission medium, the function~$u$
satisfies a diffusion evolution equation with fractional time derivative
(the precise formula will be given in~\eqref{diofhwe9ge0r8rye208tghyu:2XDAG:TRIS} below).
More explicitly,
we will perform the derivation of~\eqref{PDE}
in the space/time scale~$x := L$ and~$ t := T$
(or, more generally, at
a small spatial scale around~$ L$ and a small temporal scale around~$T$),
and the diffusive constant~$\kappa$
will be
related to some properties of the travelling profiles (namely, its slope and curvature).

{ First of all, from~\eqref{0oxc:A} and~\eqref{fo}, we see that
$$ \dot u(t)=\frac{1}{T^{2-s}}
\,
\int_0^{T^{2-s}} \dot f\left(t-\vartheta^{\frac1{2-s}}\right)\,d\vartheta
=-\frac{c}{L T^{2-s}}
\,
\int_0^{T^{2-s}} f_o'\left(\frac{x-ct+c\vartheta^{\frac1{2-s}}}{L}\right)\,d\vartheta
.$$
In particular,
\begin{equation}\label{DOT U t} 
\dot u(0)
=-\frac{c}{LT^{2-s}}
\,
\int_0^{T^{2-s}} f_o'\left(\frac{x+c\vartheta^{\frac1{2-s}}}{L}\right)\,d\vartheta
.\end{equation}
To obtain equation~\eqref{PDE} in a suitable approximation setting,
it is convenient to introduce the linear operator
$$ {\mathcal{L}}:=  \partial^s_t - {\kappa} c^s\,L^{2-s}\,\partial^2_x ,$$
for a fixed~$\kappa>0$ (to be chosen later on, see~\eqref{9qwk:9123erkasjdhg0101}).
The possibility that the linear operator~${\mathcal{L}}$ 
vanishes
is evidently equivalent to~\eqref{PDE},
and, recalling~\eqref{EG:Q3}, \eqref{EG:Q5} and~\eqref{Speed},
we have that
\begin{equation}\label{diofhwe9ge0r8rye208tghyu}
{\mathcal{L}} u(x,t)=
(2-s)\,t^{1-s}\,\dot u(0)
+ \int_0^{t^{2-s}}\ddot u (t-\vartheta^\beta) \, d\vartheta
- {\kappa} c^sL^{2-s} \,\partial^2_x u(x,t),\end{equation}
where the dependence on~$x$ is omitted whenever it creates no confusion.
Furthermore, from~\eqref{0oxc:A} and~\eqref{fo}, we see that
\begin{eqnarray*}
\partial^2_x u(x,t)&=&\frac{1}{T^{2-s}}
\,
\int_0^{T^{2-s}} \partial^2_x f\left(x,t-\vartheta^{\frac1{2-s}}\right)\,d\vartheta \\
&=&
\frac{1}{L^2 T^{2-s}}
\,
\int_0^{T^{2-s}} f_o''\left(\frac{x-ct+c\vartheta^{\frac1{2-s}}}L\right)\,d\vartheta\\
{\mbox{and }}\quad
\ddot u(t)&=&\frac{1}{T^{2-s}}
\,
\int_0^{T^{2-s}} \ddot f\left(t-\vartheta^{\frac1{2-s}}\right)\,d\vartheta\\
&=&\frac{c^2}{L^2 T^{2-s}}
\,
\int_0^{T^{2-s}} f''_o\left(\frac{x-ct+c\vartheta^{\frac1{2-s}}}L\right)\,d\vartheta.
\end{eqnarray*}
Substituting in~\eqref{diofhwe9ge0r8rye208tghyu}, 
and recalling~\eqref{DOT U t}, 
we thus find that
\begin{equation}\label{diofhwe9ge0r8rye208tghyu:2}
\begin{split} 
{\mathcal{L}} u(x,t)\,=\;&-
\frac{(2-s)\,c\,t^{1-s}}{LT^{2-s}}
\,
\int_0^{T^{2-s}} f_o'\left(\frac{x+c\vartheta^{\frac1{2-s}}}{L}\right)\,d\vartheta
\\ &\qquad+
\frac{c^2}{L^2 T^{2-s}}\int_0^{t^{2-s}}
\left[
\int_0^{T^{2-s}} f''_o\left(\frac{x-ct+c\tau^{\frac1{2-s}}+c\vartheta^{\frac1{2-s}}}L\right)\,d\vartheta
\right] \, d\tau\\ &\qquad
-\frac{{\kappa}\,c^2}{L^2}\int_0^{T^{2-s}} f_o''\left(\frac{x-ct+c\vartheta^{\frac1{2-s}}}L\right)\,d\vartheta.
\end{split}\end{equation} 
Consequently, using~\eqref{Speed} and the substitutions~$
\bar\vartheta:=T^{s-2}\vartheta$
and~$\bar\tau:=T^{s-2}\tau$,
we obtain
\begin{equation}\label{pre-diofhwe9ge0r8rye208tghyu:2theta}
\begin{split} 
{\mathcal{L}} u(x,t)\,=\;&-
\frac{(2-s)\,t^{1-s}}{T}
\,
\int_0^{1} f_o'\left(\frac{x}L+\bar\vartheta^{\frac1{2-s}}\right)\,d\bar\vartheta
\\ &\qquad+
\frac{1}{T^{s}}\int_0^{(t/T)^{2-s}}
\left[
\int_0^{1} f''_o\left(\frac{x}L-\frac{t}{T}+
\bar\vartheta^{\frac1{2-s}}+\bar\tau^{\frac1{2-s}}\right)\,d\bar\vartheta
\right] \, d\bar\tau\\ &\qquad
-\frac{{\kappa}}{T^s}\int_0^{1} 
f_o''\left(\frac{x}L-\frac{t}T+\bar\vartheta^{\frac1{2-s}}\right)\,d\bar\vartheta.
\end{split}\end{equation}
Computing \eqref{pre-diofhwe9ge0r8rye208tghyu:2theta}
at the characteristic space/time scale~$(x,t):=(L,T)$, 
we see that all the terms of~${\mathcal{L}} u(L,T)$
are of order~$T^{-s}$. More precisely,
we obtain
\begin{equation}\label{diofhwe9ge0r8rye208tghyu:2theta}
\begin{split} 
T^s\,{\mathcal{L}} u(L,T)\,=\;&-
(2-s)\,
\int_0^{1} f_o'\left(1+\bar\vartheta^{\frac1{2-s}}\right)\,d\bar\vartheta
\\&\qquad+\int_0^{1}
\left[
\int_0^{1} f''_o\left(
\bar\vartheta^{\frac1{2-s}}+\bar\tau^{\frac1{2-s}}\right)\,d\bar\vartheta
\right] \, d\bar\tau
-\kappa\int_0^{1} 
f_o''\left(\bar\vartheta^{\frac1{2-s}}\right)\,d\bar\vartheta.
\end{split}\end{equation}
Now, to detect different scales in this expression
(so to ``emphasize the quadratic part'' of the travelling wave), it
is convenient to take~$f_o$ of the form
\begin{equation}\label{PARA2}
\R\ni r\,\mapsto\, f_o(r):=a_1+a_2 r-\frac{a_3\,r^2}2+\mu\phi(r)
,\end{equation}
where~$\mu\in(0,1)$ is a small parameter, $\phi$ is a smooth function (with bounded
derivatives)
and~$a_1$, $a_2$, $a_3\in\R$.
More precisely, we will consider the case in which
\begin{equation}\label{PARA}
{\mbox{$a_2$, $a_3>0$  and~$a_3$ is sufficiently small with respect to $a_2$. }} 
\end{equation}
The case described in~\eqref{PARA} is that of a concave parabola (with sufficiently
small curvature), and the function~$f_o$ in~\eqref{PARA2} is a small perturbation
(say of size~$O(\mu)$) of such parabola, see e.g. Figure~\ref{P}.

\begin{figure}
    \centering
    \includegraphics[height=6.7cm]{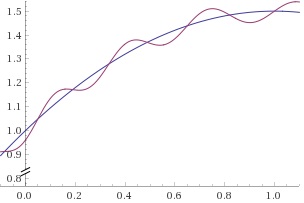}
    \caption{The function~$f_o$ in~\eqref{PARA2} and~\eqref{PARA} (a
    qualitative,
rather than quantitative, picture).}
    \label{P}
\end{figure}

In this setting, it holds that~$f_o'(r)=a_2-a_3\,r+\mu\phi'(r)$ and therefore
\begin{equation}\label{LA:LA:01o2wlLA:A12e0105j5}
\begin{split}
& f_o'\left(1+\bar\vartheta^{\frac1{2-s}}\right)
=a_2 -a_3-a_3\bar\vartheta^{\frac1{2-s}}+O(\mu)
\\
{\mbox{and }}\;\;&
f''_o\left(
\bar\vartheta^{\frac1{2-s}}+\bar\tau^{\frac1{2-s}}\right)
=-a_3+O(\mu).
\end{split}\end{equation}
Furthermore,
\begin{equation}\label{dj23ycgbvoffg43bg72346}
\int_0^{1} \bar\vartheta^{\frac1{2-s}}\,d\bar\vartheta=
\frac{2-s}{3-s}.
\end{equation}
Hence,
we insert~\eqref{LA:LA:01o2wlLA:A12e0105j5} into~\eqref{diofhwe9ge0r8rye208tghyu:2theta}, we exploit~\eqref{dj23ycgbvoffg43bg72346}
and we obtain that
\begin{equation}\label{9qwk:9123erkasjdhg01010}
\begin{split} 
T^s\,{\mathcal{L}} u(L,T)\,=\,&-
(2-s)\,
\int_0^{1} \left(
a_2 -a_3-a_3\,\bar\vartheta^{\frac1{2-s}}\right)\,d\bar\vartheta
-\int_0^{1}
a_3 \, d\bar\tau
+\kappa\int_0^{1} 
a_3\,d\bar\vartheta+O(\mu)\\
=\,& -(2-s)(a_2-a_3)+(2-s)\,a_3\cdot\frac{2-s}{3-s}-a_3+\kappa a_3
+O(\mu)\\
=\,& a_3\,\left(\kappa-(2-s)\left(\frac{a_2}{a_3}-1\right)+\frac{(2-s)^2}{3-s}-1
\right)
+O(\mu).
\end{split}\end{equation}
Therefore, the choice
\begin{equation} \label{9qwk:9123erkasjdhg0101}\kappa:=
(2-s)\left(\frac{a_2}{a_3}-1\right)-\frac{(2-s)^2}{3-s}+1 \end{equation}
makes the leading order in the right hand side
of~\eqref{9qwk:9123erkasjdhg01010} vanish. Notice also that~$\kappa>0$
in the setting of~\eqref{PARA}. This yields
that
\begin{equation}\label{diofhwe9ge0r8rye208tghyu:2XDAG:TRIS}
{\mathcal{L}} u(L,T)=O\left(\frac{\mu}{T^s}\right)=O(\delta\mu),\end{equation}
which is supposed to be a negligible quantity,
for small~$\delta$ and/or small~$\mu$.
In this approximation, we thus write~\eqref{diofhwe9ge0r8rye208tghyu:2XDAG:TRIS}
approximately as~$ {\mathcal{L}} u=0$,
thus producing
the} fractional diffusion equation in~\eqref{PDE}, as desired.

It only remains to check the claim in~\eqref{TBD}, which is the
goal of the forthcoming section.

\begin{remark}{\rm
We observe that an alternative approach to~\eqref{PDE}
consists in introducing the dimensionless variables~$\tilde t := t/T$ and~$\tilde
x := x/L$. Then, setting~$\tilde u(\tilde x,\tilde t):=u(L\tilde x,T\tilde t)$,
equation~\eqref{PDE} reduces to
$$ D^s_{\tilde t}\tilde u=\kappa\partial^2_{\tilde x}\tilde u.$$
Notice that the coefficient~$\kappa$ is expressed in
terms of constants related to some properties of the travelling profiles (e.g., according
to~\eqref{9qwk:9123erkasjdhg0101} for
a small perturbation
of a concave parabola). 

We also recall again that
such condition, as well as the derivation of~\eqref{PDE},
is obtained in~\eqref{diofhwe9ge0r8rye208tghyu:2XDAG:TRIS}
at the appropriate space/time scale~$x := L$ and~$ t := T$,
corresponding in the adimensional variables to the choice~$\tilde x=\tilde t:=1$.
More generally, one can think that
the derived
equation is valid on a small spatial scale around~$ L$ and a small temporal scale around~$
T$ (e.g. with~$T$ and~$L$ large, to make~$\delta$ small in~\eqref{DELTA}).
}\end{remark}

\section{Proof of~\eqref{TBD}}\label{STBD}

By comparing the sum with the integrals, we see that,
for each~$k\in\{1,\dots,N\}$,
$$ \int_0^{k} \frac{dx}{(x+1)^\alpha} \le \sum_{1\le j\le {k}}\frac1{j^\alpha}\le
1+\int_1^{k}\frac{dx}{x^\alpha}$$
and therefore
\begin{equation} 
\frac{(k+1)^{1-\alpha}-1}{1-\alpha}\le
\sum_{1\le j\le {k}}\frac1{j^\alpha}\le
\frac{k^{1-\alpha}-\alpha}{1-\alpha}.\end{equation}
In particular, from~\eqref{8usdhjhdgdsdf13vh:0} we have that
$$ \frac{(N+1)^{1-\alpha}-1}{1-\alpha}\le
b_N\le
\frac{N^{1-\alpha}-\alpha}{1-\alpha}.$$
Accordingly, by~\eqref{8usdhjhdgdsdf13vh},
\begin{equation}\label{DASO:0} 
\frac{\lambda_{k-1}}{L}=\frac{1}{b_N}\sum_{1\le j\le {k-1}}\frac1{j^\alpha}
\in \left[ \frac{k^{1-\alpha}-1}{N^{1-\alpha}-\alpha}\,,\;
\frac{({k-1})^{1-\alpha}-\alpha}{(N+1)^{1-\alpha}-1}
\right].\end{equation}
As a consequence,
\begin{equation}\label{DASO}
\begin{split}
& \left(\frac{k}{N}\right)^{1-\alpha}-
\frac{\lambda_{k-1}}{L} \le \left(\frac{{k}}{N}\right)^{1-\alpha}-
\frac{k^{1-\alpha}-1}{N^{1-\alpha}-\alpha}\\
&\qquad=
\left(\frac{{k}}{N}\right)^{1-\alpha} \,\left[
1-\frac{1-\frac1{{k}^{1-\alpha}} }{1-\frac{\alpha}{N^{1-\alpha}}}
\right].\end{split}
\end{equation}
Now, for large~$N$, we distinguish two cases, either~${k}\in[1,\sqrt{N}]$
or~${k}\in (\sqrt{N},N]$.
If~${k}\in[1,\sqrt{N}]$, we deduce from~\eqref{DASO} that
\begin{equation}\label{DASO2}
\left(\frac{{k}}{N}\right)^{1-\alpha}-
\frac{\lambda_{k-1}}{L} \le
\left(\frac{{k}}{N}\right)^{1-\alpha}\le
\left(\frac{\sqrt N}{N}\right)^{1-\alpha}= \frac{1}{N^{\frac{1-\alpha}2}}.
\end{equation}
If instead~${k}\in (\sqrt{N},N]$,
we deduce from~\eqref{DASO} that
\begin{equation}\label{DASO3}
\left(\frac{{k}}{N}\right)^{1-\alpha}-
\frac{\lambda_{k-1}}{L} \le
1-\frac{1-\frac1{k^{1-\alpha} }}{1-\frac{\alpha}{N^{1-\alpha}}}\le
1-\frac{1-\frac1{{N}^{\frac{1-\alpha}2}} }{1-\frac{\alpha}{N^{1-\alpha}}}.
\end{equation}
In any case, from~\eqref{DASO2} and~\eqref{DASO3}, we have that
\begin{equation}\label{DASO4}
\left(\frac{{k}}{N}\right)^{1-\alpha}-
\frac{\lambda_{k-1}}{L} \le\frac{1}{N^{\frac{1-\alpha}2}}+
\left|
1-\frac{1-\frac1{{N^{\frac{1-\alpha}2}}} }{1-\frac{\alpha}{N^{1-\alpha}}}\right|.\end{equation}
In addition, from~\eqref{DASO:0} we also have that
\begin{equation}\label{BASO1}
\begin{split}&
\frac{\lambda_{k-1}}{L} -
\left(\frac{{k}}{N}\right)^{1-\alpha}\le
\frac{({k-1})^{1-\alpha}-\alpha}{(N+1)^{1-\alpha}-1}-
\left(\frac{{k}}{N}\right)^{1-\alpha}\\
&\qquad\le
\left(\frac{{k}}{N}\right)^{1-\alpha}
\left[ \frac{\left(\frac{k-1}{k}\right)^{1-\alpha}-\frac\alpha{k^{1-\alpha}}}{
\left(\frac{N+1}{N}\right)^{1-\alpha}-\frac1{N^{1-\alpha}}}
-1
\right].
\end{split}
\end{equation}
Now, if~$k\in[1,\sqrt{N}]$,
we infer from~\eqref{BASO1} that
\begin{equation}\label{BASO2}
\begin{split}&
\frac{\lambda_{k-1}}{L} -
\left(\frac{{k}}{N}\right)^{1-\alpha}\le 2\,
\left(\frac{{k}}{N}\right)^{1-\alpha}\le 2\,
\left(\frac{\sqrt{N}}{N}\right)^{1-\alpha}=
\frac{2}{N^{\frac{1-\alpha}2}}.
\end{split}
\end{equation}
If instead~${k}\in (\sqrt{N},N]$, we deduce from~\eqref{BASO1} that
\begin{eqnarray*} &&
\frac{\lambda_{k-1}}{L} -
\left(\frac{{k}}{N}\right)^{1-\alpha}\le 
\frac{\left(\frac{k-1}{k}\right)^{1-\alpha}-\frac\alpha{k^{1-\alpha}}}{
\left(\frac{N+1}{N}\right)^{1-\alpha}-\frac1{N^{1-\alpha}}}
-1\le
\frac{1-\frac\alpha{N^{1-\alpha}}}{1-\frac1{N^{1-\alpha}}}
-1.
\end{eqnarray*}
In any case, recalling~\eqref{BASO2}, we have that for every~$k\in\{1,\dots,N\}$
it holds that
$$ \frac{\lambda_{k-1}}{L} -
\left(\frac{{k}}{N}\right)^{1-\alpha}\le \frac{2}{N^{\frac{1-\alpha}2}}+
\left|
\frac{1-\frac\alpha{N^{1-\alpha}}}{1-\frac1{N^{1-\alpha}}}
-1\right|.
$$
Hence, in view of~\eqref{uaso:n} and~\eqref{DASO4}, we have that
\[
\frac{|\eta_{k,N}|}{ L}={
\left|\left(\frac{k}{N}\right)^{1-\alpha}-\frac{\lambda_{k-1}}L\right| }
\le\frac{3}{N^{\frac{1-\alpha}2}}+
\left|
1-\frac{1-\frac1{{N^{\frac{1-\alpha}2}}} }{1-\frac{\alpha}{N^{1-\alpha}}}\right|+
\left|
\frac{1-\frac\alpha{N^{1-\alpha}}}{1-\frac1{N^{1-\alpha}}}
-1\right|,
\]
and this plainly implies~\eqref{TBD}, as desired.

\section{Conclusions}\label{STBD:2}

Purkinje cells seem to exhibit two special features: 
\begin{enumerate}
\item[(i)] \hskip0.5cm on the one hand, their malfunction is related, among the others, to
abrupt movements, tremors and lack of coordination;
\item[(ii)] \hskip0.5cm on the other hand, recent experiments have shown the evidence
of time fractional diffusion arising in Purkinje cells.
\end{enumerate}
Also, in the mathematical theory of evolution equations, typically two regularity regimes
arise:
\begin{enumerate}
\item[(i)'] \hskip0.5cm on the one hand, hyperbolic equations typically
present shocks and
irregular solutions;
\item[(ii)'] \hskip0.5cm on the other hand, parabolic equations are endowed with
good regularity theories with respect to the initial data.
\end{enumerate}
It is quite tempting to relate the biological phenomena in (i)--(ii)
with the mathematical treats in (i)'-(ii)', respectively.
In this paper, we provide a mathematical setting to show how
highly ramified media affect the propagation of travelling waves,
by producing a superposition of delayed signals which in turn may
be related to time fractional diffusion.

In our computation, we derive the time fractional heat equation 
as the effect of the superposition of delayed classical heat equations,
measured at a spatial scale in proximity of the end of the transmission network
and at a time close to the characteristic diffusion time. 
The delay in the fundamental heat equation is produced
by a network, in which the length of each branch is appropriately chosen to produce
a suitable superposition effect.

Our expansion is
calculated on a specific travelling wave, taken as a small perturbation of
a concave parabola, and the diffusion coefficient
depends on the slope and on the curvature of such
travelling profile.

The remainders of our formula are explicitly stated and the asymptotics are discussed
in details. To emphasize the role
of the scales at which the main equation is attained (up to remainders),
one can also rewrite the dimensional
time fractional heat equation in terms of adimensional space and time variables.

In our mathematical framework, we show that
highly ramified media may provide a regularizing effect
on the leading equations of signal transmission.
It is of course intriguing
to relate the ramification of these underlying mathematical media
with the structure of the
Purkinje neuron's dendritic arbor. 

We also recall that
usually fractional diffusion is discussed mostly in relation
to its classical parabolic analogue (for instance,
it is commonly viewed that ``the appearance of fractional equations
is very appealing due to their proximity to the analogous standard equations'',
see page~5 in~\cite{YY77}). In this sense,
from the theoretical point of view, 
our approach seems to be rather different from the existing literature and
the new link that we propose
between
fractional diffusion and hyperbolic (rather than parabolic) equations
may lead to stimulating mathematical
considerations from a different perspective. 

\section*{Acknowledgements}

It is a pleasure to thank Elena Saftenku and
Fidel Santamaria for very interesting discussions
on neural transmissions. {We also thank the Referees for their very
valuable comments.}

This work has been supported by the Australian Research Council 
Discovery Project ``N.E.W.
Nonlocal Equations at Work''.

\vfill


\begin{thebibliography}{99}

\bibitem{ABAR}
{\sc N. Abatangelo, E. Valdinoci}: Getting
acquainted with the fractional Laplacian.
To appear in 
Springer INdAM Ser.

\bibitem{allen}
{\sc M. Allen}:
A Nondivergence Parabolic Problem with a Fractional Time derivative.
{\em Differential Integral Equations} {\bf 31}, no. 3-4
(2018), 215--230. 

\bibitem{allen-et}
{\sc M. Allen, L. Caffarelli, A. Vasseur}:
A parabolic problem with a fractional time derivative.
{\em Arch. Ration. Mech. Anal.} {\bf 221}, no. 2 (2016), 603--630.

\bibitem{ANAS}
{\sc T. J. Anastasio}:
Nonuniformity in the linear network model of the oculomotor
integrator produces approximately fractional-order dynamics and more
realistic neuron behavior.
{\em Biol Cybern.} {\bf 79} (1998), 377--391.

\bibitem{RR1}
{\sc R. Appali, U. van Rienen, T. Heimburg}:
A comparison of the Hodgkin-Huxley model the
soliton theory for the action potential in nerves. 
{\em Adv. Planar Lipid Bilayers Liposomes} {\bf 16} (2012), 275--298.

{
\bibitem{EQUI}
{\sc R. Bagley}:
On the equivalence of the Riemann-Liouville and the Caputo fractional order
derivatives in modeling of linear viscoelastic materials.
{\em Fract. Calc. Appl. Anal.} {\bf 10} (2007), no. 2, 123--126. 
}

\bibitem{COST}
{\sc C. A. Balanis}:
{\em Advanced engineering electromagnetics}.
John Wiley \& Sons, 2012.

\bibitem{ben}
{\sc D. ben-Avraham, S. Havlin}:
{\em Diffusion and reactions in fractals and disordered systems}.
Cambridge University Press, Cambridge, 2000. 

\bibitem{blanco}
{\sc A. Blanco, R. Moyano, J. Vivo,
R. Flores-Acu\~{n}a, A. Molina, C. Blanco, J. G. Monterde}:
Purkinje cell apoptosis in arabian horses with cerebellar abiotrophy.
{\em J. Vet. Med. Physiol. Pathol. Clin. Med.} {\bf 53}, no. 6 (2006), 286--287.

\bibitem{bucur}
{\sc C. Bucur}:
{\em Local density of Caputo-stationary functions in the space of
smooth functions}.
{\em ESAIM Control Optim. Calc. Var.} {\bf 23}, no. 4 (2017), 1361--1380. 

\bibitem{bucur2}
{\sc C. Bucur, E. Valdinoci}:
{\em Nonlocal diffusion and applications}.
Lecture Notes of the Unione Matematica Italiana 20.
Springer, Bologna, 2016.

\bibitem{caputo}
{\sc M. Caputo}: 
Linear model of dissipation whose $Q$ is 
almost frequency independent-II.
{\em Geophys. J. R. Astron. Soc.} {\bf 13}, no. 5 (1967), 529--539.

\bibitem{COO}
{\sc S. Coombes}:
Neural fields.
Scholarpedia {\bf 1}, no. 6 (2006), 1373.

\bibitem{zuazua}
{\sc R. D\'ager, E. Zuazua}:
{\em Wave propagation, observation and control in $1$-d
flexible multi-structures}. 
Math\'ematiques \& Applications, 50. 
Springer-Verlag, Berlin, 2006. 

\bibitem{kai}
{\sc K. Diethelm}:
{\em The analysis of fractional differential equations. 
An application-oriented exposition using differential operators of Caputo type}.
Lecture Notes in Mathematics. Springer, Berlin, 2004.

\bibitem{JJ}
{\sc S. Dipierro, O. Savin, E. Valdinoci}:
All functions are locally $ s $-harmonic up to a small error.
{\em J. Eur. Math. Soc. (JEMS)} 
{\bf 19} (2017), no. 4, 957--966.

\bibitem{JJ2}
{\sc S. Dipierro, O. Savin, E. Valdinoci}:
Local approximation of arbitrary functions by solutions of nonlocal equations.
ArXiv 1609.04438 (2016).

{
\bibitem{DUHU}
{\sc M. Du, Z. Wang, H. Hu}.
Measuring memory with the order of fractional derivative.
{\em Sci. Rep.} {\bf 3} (2013), 3431.
}

\bibitem{EL}
{\sc A. El Hady, B. B. Machta}:
Mechanical surface waves accompany
action potential propagation. {\em Nature Comm.} {\bf 6} (2015),
6697 EP.

\bibitem{ERME01}
{\sc G. B. Ermentrout, D. Kleinfeld}:
Traveling electrical waves in cortex: Insights from phase dynamics 
and speculation on a computational role. {\em Neuron} 29 (2001), 33--44.

\bibitem{ERME93}
{\sc G. B. Ermentrout, J. B. McLeod}:
Existence and uniqueness of travelling waves for a neural network. 
{\em Proc. Royal Soc. Edinburgh} 123A (1993), 461--478.

\bibitem{evans}
{\sc L. C. Evans}:
{\em Partial differential equations}. Graduate Studies in
Mathematics, 19. American Mathematical Society, Providence, RI, 1998.

\bibitem{fiala}
{\sc J. C. Fiala, K. M. Harris}:
Dendrite structure. In {\sc
G. Stuart, S. Nelson, and M. H\"ausser, eds.}:
{\em Dendrites}. Oxford Scholarship Online.
Oxford University
Press, Oxford, 1999.

\bibitem{GGG}
{\sc A. Gonzalez-Perez, L. D. Mosgaard, 
R. Budvytyte, E. Villagran-Vargas, A. D. Jackson, T. Heimburg}: 
Solitary electromechanical pulses in lobster neurons, {\em
Bioph. Chem.} {\bf 216} (2016), 51--59.

\bibitem{RR2}
{\sc T. Heimburg, A. D. Jackson}:
On soliton propagation in biomembranes and nerves. {\em Proc. Nat.
Acad. Sci.} {\bf 102}, no. 28 (2005), 9790--9795.

\bibitem{HHH}
{\sc A. L. Hodgkin, A. F. Huxley}:
A quantitative description of membrane current
and its application to conduction and excitation in nerve.
{\em J. Physiology} {\bf 117}, no.~4 (1952), 500--544.

{
\bibitem{IONE}
{\sc C. Ionescu, A. Lopes, D. Copot, J. A. T. Machado, J. H. T. Bates}:
The role of fractional calculus in modelling biological phenomena:
a review. 
{\em Comm. Nonlin.
Science Numer. Simul.} 
{\bf 51} (2017), 141--159.
}

\bibitem{RR3}
{\sc V. G. Ivancevic, T. T. Ivancevic}: 
{\em Quantum Neural Computation}. Springer, Netherlands, 2010.

\bibitem{KARA}
{\sc J. J. Karakash}:
{\em Transmission Lines and Filter Networks}. Macmillan,  New York, 1950.

\bibitem{kim}
{\sc I. Kim, K.-H. Kim, S. Lim}:
An $L_q (L_p )$-theory for the time fractional
evolution equations with variable coefficients.
{\em Adv. Math.} {\bf 306} (2017), 123--176. 

\bibitem{LARKUM}
{\sc M. E.
Larkum, S. Watanabe, T. Nakamura, N. Lasser-Ross, W. N. Ross}:
Synaptically activated Ca2+ waves in
layer 2/3 and layer 5 rat neocortical pyramidal neurons.
{\em J. Physiol.} {\bf 549}
(2003), 471--488.

\bibitem{RR4}
{\sc B. Lautrup, R. Appali, A. D. Jackson, T. Heimburg}:
The stability of solitons in biomembranes
and nerves. {\em Eur. Phys. J.} {\bf 34}, no. 57 (2011), 1--9.

\bibitem{Magin}
{\sc R. L. Magin}:
Fractional calculus models of complex dynamics in biological tissues.
{\em Computers Math. Appl.} {\bf 59} (2010), 1586--1593.

\bibitem{Fidel2}
{\sc T. Marinov, F. Santamaria}:
Modeling the effects of anomalous diffusion on
synaptic plasticity.
{\em BMC Neuroscience} {\bf 14}, Suppl. 1 (2013),
P343.

\bibitem{MARIN}
{\sc T. Marinov, F. Santamaria}:
Computational modeling of diffusion in the cerebellum.
{\em Prog. Mol. Biol. Transl. Sci.} {\bf 123} (2014), 169--89.

\bibitem{MF}
{\sc I. A. Mavroudis, D. F. Fotiou, L. F. Adipepe, M. G. Manani, 
S. D. Njau, D. Psaroulis, V. G. Costa, S. J. Baloyannis}:
Morphological changes of the human purkinje cells and 
deposition of neuritic plaques and neurofibrillary tangles 
on the cerebellar cortex of Alzheimer's disease. 
{\em Amer. J. Alzheimer's Dise. Other Dem.} {\bf 25}, no. 7 (2010), 585--591.

\bibitem{YY77}
{\sc R. Metzler, J. Klafter}:
The random walk’s guide to anomalous diffusion: 
A fractional dynamics approach.
{\em Phys. Rep.} {\bf 339}, no.~1 (2000), 1--77. 

\bibitem{MIRA} {\sc W. L. Miranker}:
A neural network wave formalism.
Adv. Appl. Math., {\bf 37} (2006) 19--30.

\bibitem{NACA}
{\sc T. Nakamura, N. Lasser-Ross, K. Nakamura, W. N. Ross}:
Spatial segregation and interaction of calcium
signalling mechanisms in rat hippocampal CA1 pyramidal neurons.
{\em J. Physiol.} {\bf 543}
(2002),
465--480.

\bibitem{SSS}
{\sc S. A. Neymotin, R. A. McDougal,
M. A. Sherif, C. P. Fall, M. L. Hines, W. W. Lytton}:
Neuronal Calcium wave propagation varies with changes in 
endoplasmic reticulum parameters: a computer model.
{\em Neural Comp.} {\bf 27}, no.~4 (2015),
898--924.

\bibitem{2002}
{\sc E. A. Nimchinsky, B. L. Sabatini, K. Svoboda}: Structure and function 
of dendritic spines. {\em Annu. Rev. Physiol.} {\bf 64} (2002),
313--353.

\bibitem{ERME01b}
{\sc D. J. Pinto, G. B. Ermentrout}:
Spatially structured activity in synaptically
coupled neuronal networks: I. Travelling fronts and pulses. {\em 
SIAM J.
Applied Math.}, {\bf 62} (2001), 206--225.

\bibitem{RIGA}
{\sc G. G. Rigatos}:
{\em Advanced Models of Neural Networks}.
Nonlinear Dynamics and Stochasticity in Biological Neurons.
Springer-Verlag, Berlin, 2015.

\bibitem{ROSS}
{\sc W. N. Ross}:
Understanding calcium waves and sparks in central neurons.
{\em 
Nature Rev. Neurosci.} {\bf 13} (2002), 157--168.

\bibitem{EE}
{\sc E. \`E. Saftenku}:
Modeling of slow glutamate diffusion and AMPA receptor activation
in the cerebellar glomerulus.
{\em J. Theor. Biology} {\bf 234} (2005), 363--382.


\bibitem{Fidel}
{\sc F. Santamaria, S. Wils, E. De Schutter, G. J. Augustine}:
Anomalous Diffusion in Purkinje
Cell Dendrites Caused by Spines.
{\em Neuron} {\bf 52} (2006), 635--648.

\bibitem{Fidel3}
{\sc F. Santamaria, S. Wils, E. De Schutter, G. J. Augustine}:
The diffusional
properties of dendrites depend on the density of dendritic spines.
{\em Europ. J. Neuroscience} {\bf 34}, no.~4 (2011), 561--568.

\bibitem{BIND}
{\sc M. J. Saxton}: Anomalous diffusion due to binding: a Monte
Carlo study. {\em Biophys. J.} {\bf 70} (1996), 1250--1262.

\bibitem{Schweidler}
{\sc E. R. von Schweidler}:
Studien \"uber die Anomalien im Verhalten der Dielectrika.
{\em Ann. Phys.} {\bf 24}, (1907), 711--770.

\bibitem{Thorson}
{\sc J. Thorson, M. Biederman-Thorson}:
Distributed relaxation processes
in sensory adaptation: spatial nonuniformity in receptors can explain both the
curious dynamics and logarithmic statics of adaptation.
{\em Science} 
{\bf 183}, no.~4121 (1974), 161--172.

\bibitem{44}
{\sc J. Trommersh\"auser, J. Marienhagen, A. Zippelius}:
Stochastic Model of Central Synapses: Slow Diffusion of 
Transmitter Interacting with Spatially Distributed Receptors and Transporters.
{\em J. Theor. Biology} {\bf 198} (1999), 101--120.

\bibitem{WIKI} {\sc Wikipedia}:
Drawing of Purkinje cells (A) and granule cells (B)
from pigeon cerebellum by Santiago Ram\'on y Cajal, 1899;
Instituto Cajal, Madrid, Spain.
File:PurkinjeCell.jpg
\\
{\tt https://en.wikipedia.org/wiki/Purkinje\_cell\#/media/File:PurkinjeCell.jpg}

\bibitem{zacher}
{\sc R. Zacher}:
Maximal regularity of type $L_p$ for abstract parabolic Volterra equations. 
{\em J. Evol. Equ.} {\bf 5}, no. 1 (2005), 79--103. 

\bibitem{zacher2}
{\sc R. Zacher}: A De Giorgi-Nash type theorem for time
fractional diffusion equations. {\em Math.
Ann.} {\bf 356}, no. 1 (2013), 99--146.

\end{thebibliography}
\end{document}